\documentclass[11pt,leqno,fleqn]{article}
\newenvironment{ISItext}
{\normalsize\rm\setlength{\parindent}{1cm}\setlength{\parskip}{0pt}}
{\vskip 12pt}

\newenvironment{ISIresume}
{\normalsize\rm\setlength{\parindent}{1cm}\setlength{\parskip}{0pt}\it}
{\vskip 12pt}

\begin{document}
\newcommand{\ISItitle}[1]{\vskip 0pt\setlength{\parindent}{0cm}\Large\textbf{#1}\vskip 12pt}
\newcommand{\ISIsubtitleA}[1]{\normalsize\rm\setlength{\parindent}{0cm}\textbf{#1}\vskip 12pt}
\newcommand{\ISIsubtitleB}[1]{\normalsize\rm\setlength{\parindent}{0cm}\textbf{#1}\vskip 12pt}
\newcommand{\ISIsubtitleFig}[1]{\normalsize\rm\setlength{\parindent}{0cm}
\textbf{\textit{#1}}\vskip 12pt}
\newcommand{\ISIauthname}[1]{\normalsize\rm\setlength{\parindent}{0cm}#1 \\}
\newcommand{\ISIauthaddr}[1]{\normalsize\rm\setlength{\parindent}{0cm}\it #1 \vskip 12pt}

\ISItitle{First Passage Times and Breakthrough Curves Associated with Interfacial Phenomena}

\ISIauthname{Waymire, Edward C} 
\ISIauthaddr{Oregon State University, Mathematics\\
Kidder Hall\\
Corvallis 97331, USA\\ 
E-mail: waymire@math.oregonstate.edu}

\ISIauthname{Appuhamillage, Thilanka A } 
\ISIauthaddr{Oregon State University, Mathematics\\
Kidder Hall\\
Corvallis 97331, USA\\ 
E-mail: ireshara@math.oregonstate.edu}

\ISIauthname{Bokil, Vrushali A }
\ISIauthaddr{Oregon State University, Mathematics\\
Kidder Hall\\
Corvallis 97331, USA\\ 
E-mail: bokilv@math.oregonstate.edu}


\ISIauthname{Thomann, Enrique A }
\ISIauthaddr{Oregon State University, Mathematics\\
Kidder Hall\\
Corvallis 97331, USA\\ 
E-mail: thomann@math.oregonstate.edu}

\ISIauthname{Wood, Brian D  }
\ISIauthaddr{Oregon State University, Mathematics\\
Kidder Hall\\
Corvallis 97331, USA\\ 
E-mail: wood@math.oregonstate.edu}

\begin{ISItext}
The topic addressed in this paper was initially motivated by 
questions resulting from recent
laboratory experiments designed to empirically test
and understand advection-dispersion in the presence of sharp interfaces; e.g., 
experiments by
\cite{kuo_1999}, 
\cite{hoteit} 
\cite{Berkowitz09}. 
Such laboratory experiments have 
been rather sophisticated in the use of layers of sands and/or
glass beads of different granularities  and modern measurement technology.  As a result they have uncovered
a convincing empirical foundation for some interesting and unexpected phenomena that had
escaped prior theoretical notice and explanation.  
To this end it is natural to inquire about the effect of an
interface on the stochastic particle motion of  immersed solutes.    From a general 
mathematical point of view an {\it interface} is defined by a hypersurface across which the dispersion coefficient is discontinuous. 
As is well-known for the case  of dilute suspensions
in a homogeneous medium (e.g., water),  perhaps flowing at a rate $v$,  the particle motion is that of  a  Brownian motion with a constant diffusion
coefficient $D > 0$ and drift $v$.   For simple one-dimensional flow across an interface,  a localized point interface results in 
a skewness effect that explains much of the empirically observed results noted above; see 
\cite{Ramirez06}, \cite{Ramirez08},
\cite{App09a}, \cite{App09b}, \cite{Ramirez11a}. 

\vskip 5in

\noindent {\bf Problem}  As an illustration of empirical findings, suppose that a dilute solute is injected at a point $L$ units to the left of 
an interface at the origin and retrieved at a point $L$ units to the right of the interface.  Let $D^-$ denote
the (constant) dispersion coefficient to the left of the origin and $D^+$ that to the right, with
say $D^- < D^+$.   Conversely,
suppose the solute is injected at a point $L$ units to the right of the interface and retrieved at a point $L$
units to the left.  {\it Which will be retrieved first ? } 

For positive parameters $D^+, D^-$, consider a piecewise constant dispersion coefficient with interface
at $x = 0$ given by
$$D(x) = D^-1_{(-\infty,0)}(x) + D^+1_{[0,\infty)}(x), \quad x\in {\bf R}.$$

\vskip .15in
\noindent{\bf Theorem 1}
 Let $D^+, D^-$ be arbitrary positive numbers, with say  $D^- < D^+$.    
Define $Y^{(\alpha^*)}_t = s(B_t^{(\alpha^*)}), t\ge 0$, where $B^{(\alpha^*)}$ is skew Brownian motion with transmission
parameter $\alpha^* =  {D^+\over D^+ + D^-}$, and
 $s(x) = \sqrt{D^+}x1_{[0,\infty)}(x) + \sqrt{D^-}x1_{(-\infty,0]}(x), x\in {\bf R}$.
Let $T_y = \inf\{t\ge 0:  Y^{(\alpha^*)}_t  = y\}.$  Then,

\vskip .05in

\noindent (a) For smooth  initial data  $c_0$,   \    $c(t,y) = E_yc_0(Y^{(\alpha)}_t ), \quad t\ge 0$, solves 
$${\partial c\over\partial t} = {1\over 2}{\partial\over\partial y}(D(y){\partial c\over\partial y}),\quad D^+{\partial c(t,0^+)\over\partial y} = D^-{\partial c(t,0^-)\over\partial y}.$$

\vskip .05in

\noindent (b)  For $y > 0$,  $\quad P_{-y}(T_y > t) \le {\sqrt{D^-}\over\sqrt{D^+}} P_{y}(T_{-y} > t) <  P_{y}(T_{-y} > t) , \quad t\ge 0.$

\vskip .05in

\vskip .15in
\noindent{\bf Remark}  {\it This basic result was obtained \cite{App09b} in terms of first passage times, however
the factor $\sqrt{D^-}/\sqrt{D^+}$ was not included in the statement of the result there.   Related phenomena and  results on dispersion 
in this context are also given in \cite{Ramirez06}, \cite{Ramirez08}, \cite{App09a}, \cite{Ramirez11a}.   In addition, a formula for the first passage time distribution for skew Brownian
motion was recently obtained in \cite{AppShel}.  In principle, the identification of stochastic particle motions can have computational advantages.
Results pertaining to Monte-Carlo simulations of skew diffusions are described in 
\cite{Lejay08} and references therein.}

\vskip .15in

As illustrated by the examples below, the role of interfacial phenomena is of much broader
interest than suggested by  advection-dispersion experiments.  However the specific 
nature of the interface can vary, depending on the specific phenomena.    We briefly describe three distinct classes
of examples  of phenomena from the biological/ecological sciences in which interfaces naturally occur.

\vskip .05in

\noindent{\bf Example 1} ({\it\bf Coastal  Upwelling and Fisheries})  Up-wellings,  the movement of deep nutrient 
rich waters to the sun-lit ocean surface,  occur in roughly
one percent of the ocean but are responsible for nearly fifty-percent of the worlds fishing industry.
The up-welling along the  Malvinas current that occurs off of the coast of Argentina is unusual
in that it is the result of  a very sharp break in the shelf, rather than being driven by winds.
 The equation for the free surface $\eta$
as a function of spatial variables $(x,y)$ is of the form
$$
{\partial\eta\over\partial y} =  -{r\over f}\left({\partial h\over\partial x}\right)^{-1}{\partial^2\eta\over\partial x^2},$$
where $r > 0$ and $f < 0$ in the  southern hemisphere, and $h(x)$ is the depth of the ocean at a distance $x$ from 
the shore.  In particular, the sharp break in the shelf makes $h^\prime(x)$ a piecewise constant function
with  positive values $H^+, H^-$.  The location of the interface coincides with the  distance to the shelf-break.
If the spatial variable $y > 0$ is viewed as a \lq\lq time\rq\rq  parameter, then this is a skew-diffusion equation,
however the physics imply continuity of the derivatives $\partial\eta/\partial y$ at the interface; see \cite{Matano} and
references therein.

\vskip 5in
\noindent{\bf Example 2} ({\it\bf Fender's Blue Butterfly}) The Fenders Blue is an endangered species of butterfly found
in the pacific northwestern United States.  The primary habitat patch is Kinkaid's Lupin flower.   Quoting 
\cite{Schultz},  {\it \lq\lq Given past research on the Fender's blue,  and the potential to investigate response to
patch boundaries, we ask two central questions.  First,  how do organisms respond to habitat edges? Second,
what are the  implications of this behavior  for residence times?\rq\rq}  Sufficiently long residence (occupation) times
in Lupin patches  are required for pollination, eggs, larvae and ultimate sustainability of the population. 
Empirical evidence points  to a skewness in random walk models for butterfly movement at the path boundaries.

\vskip .05in
\noindent{\bf Example 3} ({\it\bf Sustainability on a River Network}) The movement of larvae in a river system is often modeled by
advective-dispersion equations  in which the rates are determined by hydrologic/geomorphologic relationships in the form
of the  so-called {\it Horton laws}.   In general river networks are modeled as directed binary tree graphs and each junction
may be viewed as an interface.   Conservation of mass leads to continuity of  flux of larvae across each stream junction as
the  appropriate interface condition.   Problems on sustainability in this context are generally formulated in terms of 
network size and characteristics relative to the production of larvae sufficient to prevent permanent downstream removal
at low population sizes; see 
 \cite{Ramirez11b} for recent results in the case of a river network.

 \vskip .05in
 
The following theorem provides a useful summary of the interplay between diffusion coefficients and broader classes of possible 
interfacial conditions illustrated by these examples.    The  proof follows by  a straightforward application of the It\^o-Tanaka formula.

\vskip .15in
\noindent{\bf Theorem 2} Let $D^+, D^-$ be arbitrary positive numbers and let  $0 <\alpha,  \lambda < 1$.    
Define $Y^{(\alpha)}_t = s(B_t^{(\alpha})), t\ge 0$, where $B^{(\alpha)}$ is skew Brownian motion with transmission
parameter $\alpha$ and $s(x) = \sqrt{D^+}x1_{[0,\infty)}(x) + \sqrt{D^-}x1_{(-\infty,0]}(x), \ x\in {\bf R}$.
Then 
$$M_t = f(Y^{(\alpha)}_t) - {1\over 2}\int_0^t D(Y^{(\alpha)}_u)f^{\prime\prime}(Y_u)du, \quad t\ge 0,$$
is a martingale for all $f\in {\cal D}_\lambda = \{f\in C^2({\bf R}\backslash \{0\})\cap C({\bf R}) :  \lambda f^\prime(0^+) 
= (1-\lambda)f^\prime(0^-)\}$ if and only if 
$$\alpha = \alpha^*(\lambda)  = { \lambda\sqrt{D^-}\over \lambda\sqrt{D^-} + (1-\lambda)\sqrt{D^+}}.$$

\vskip .15in

\noindent{\bf Remark} {\it This theorem is a generalization of the results obtained by 
\cite{Ramirez06} and \cite{App09a}
for the case of advection-dispersion problems across an interface described at the outset, where the parameter $\lambda = {D^+\over D^+ + D^-}$ and
$\alpha^* = {\sqrt{D^+}\over\sqrt{D^+} + \sqrt{D^-}}.$ }

\vskip .15in

\noindent{\bf Definition} With the choice of $\alpha^* \equiv\alpha^*(\lambda)$ given by Theorem 2, we
 refer to the process $Y^{\alpha^*}$ as the {\it physical diffusion} corresponding to the 
dispersion coefficients $D^+, D^-$ and interface parameter $\lambda$. 

\vskip .15in

Observe that in the application to the coastal up-welling problem one obtains
$$\alpha^* = {\sqrt{D^-}\over\sqrt{D^+} + \sqrt{D^-}}.$$   The physical diffusion for this
example may be checked to
coincide with the Stoock-Varadahn martingale in this case; see \cite{SV} for the definition
of the corresponding martingale
problem.  Note that the answer to the first passage time problem will be exactly opposite
to that obtained for advection-dispersion experiments under this model.

\vskip 5in

We conclude with a result to show that the issue raised in Example 2 relating interfacial conditions
to residence times is indeed a sensitive problem.

\vskip .15in
\noindent{\bf  Theorem 3} Let $Y^{\alpha^*}$  denote the physical diffusion for the dispersion 
coefficients $D^+, D^-$ and interface parameter $\lambda$. Define modified occupation time
processes by
$$\tilde{\Gamma}^+(t) = \int_0^t 1[Y^{(\alpha^*)}_s > 0]ds, \quad t\ge 0.$$
Similarly let $\tilde{\Gamma}^-(t) = t- \tilde{\Gamma}^+(t), t\ge 0.$  Then,
$$\tilde{\Gamma}^+(t) > \tilde{\Gamma}^-(t) \ \forall t>0\quad {\iff} \lambda > {\sqrt{D^+}\over\sqrt{D^+} +\sqrt{D^-}},$$
with equality  when $\lambda = {\sqrt{D^+}\over\sqrt{D^+} +\sqrt{D^-}}$.
\vskip .15in
\noindent{\it Proof}
Let $\lambda(x) = 2\lambda1_{[0,\infty)}(x) + 2(1-\lambda)1_{(-\infty,0)}(x)$, and define $\rho(x) = D(x)/\lambda(x)$.  Consider the time
change $\tau_\rho(t,\omega)$ defined by
$$\int_0^{\tau_\rho(t)} {1\over \rho(B_s)} = t, \quad t\ge 0.$$
Define $S_\rho$ by $B(t,S_\rho(\omega)) = B(\tau_\rho(t,\omega),\omega)
= Z(t,\omega)$.   Then the process $Z$ is a diffusion with zero drift and 
diffusion coefficient $D^\pm\over \lambda^2(x)$ with interface parameter $\lambda = 1/2$.
Now observe that
\begin{eqnarray*}
\Gamma^+_{Y^{(\alpha^*)}}(t) &:=& \int_0^t1[Y^{(\alpha^*)}_s>0]d<Y^{(\alpha^*)}>_s\\
&=& \int_0^t1[Z_s>0]4\lambda^2d<Z>_s\\
&=& 4\lambda^2\int_0^t1[B(\tau_\rho(s)>0]{D^+\over 4\lambda^2}ds\\
&=& 4\lambda^2\Gamma^+_B(\tau_\rho(t)).
\end{eqnarray*}
Similarly $\Gamma^{-}_{Y^{(\alpha^*)}}(t) = 4(1-\lambda)^2\Gamma^-_B(\tau_\rho(t))$.  Thus $(1-\lambda)^2\Gamma^{+}_{Y^{(\alpha^*)}}(t)
= \lambda^2\Gamma^{-}_{Y^{(\alpha^*)}}(t).$  Now observe that
$$\tilde\Gamma^+_{Y^{(\alpha^*)}}(t) = \int_0^t1[Z_s>0]ds
= {4\lambda^2\over D^+}\Gamma_B^+(\tau_\rho(t)),$$
and similarly for $\tilde\Gamma^+_{Y^{(\alpha^*)}}(t)$, to arrive at
$${D^-\over (1-\lambda)^2}\tilde{\Gamma}^-(t)  = {D^+\over \lambda^2}\tilde\Gamma^+(t).$$
The assertion now follows.
\rightline{QED}

\vskip .2in

It is interesting to note that under the mass conservation
 interface parameter $\lambda = D^+/(D^+ + D^-)$, the particle will reside longer in the region
 with the faster dispersion rate.  While this is to be expected for physical experiments of dispersion
 in porous media of the type described above,  it shows that the conservative interface condition (defined
 by this choice of $\lambda$)  is likely not appropriate for models of animal movement !

\vskip 5in

\noindent{\bf Remark}  {\it A related phenomena in terms of a \lq\lq modified local time\rq\rq is described
in \cite{App09b}.  The modification, denoted with  the $\tilde{}$, refers to an integration with respect to Lebesgue measure in place of 
quadratic variation in the usual mathematical definition of local time and quadratic variation; e.g., see
\cite{RY}.
 In general, the treatment of dispersion in the presence of interfaces suggests that
the physical/biological theories are, to the  extent possible, 
 naturally based on a modification of local times and occupation times
in which integration with respect to quadratic variation is replaced by integration with respect to 
Lebesgue measure.  In fact it is shown in \cite{App09b} that this naturally leads to a stochastic 
determination of the physical transmission parameter $\alpha^*$ in terms of a continuity condition
on the modified local time of the stochastic particle motion.  This is a probabilistic condition at
the  particle scale that may be viewed as an alternative to the usual macro-scale  pde condition of continuity of flux in particle  concentrations.}

\vskip .45in

\ISIsubtitleB{ACKNOWLEDGMENTS}
The results of this paper appear in the extended abstracts of the 58th World Congress of the  International Statistics Institute, Dublin, Ireland, August 21--26, 2011.
  This is partially supported by a grant from the National Science Foundation, proposal number EAR-0724865.

\end{ISItext}








\ISIsubtitleB{ABSTRACT}

\begin{ISIresume}
Advection and dispersion in highly heterogeneous environments
involving interfacial discontinuities in the corresponding drift and dispersion rates
are described through disparate examples from the physical and biological sciences.  A mathematical framework is formulated to address
specific empirical phenomena involving first passage time and occupation time functionals observed in relation to the  interfacial parameters.

\end{ISIresume}


\begin{thebibliography}{plain}
\begin{small}
\bibitem{App09a}
 Appuhamillage, T., V. Bokil, E. Thomann, E. Waymire, B. Wood (2011)
  Occupation and Local Times for Skew Brownian Motion with Applications to Dispersion Across an Interface, 
  Annals of Appld. Probab,   {\bf  21}(1) 183--214. [Correction: Ann. Appld. Probab., 
  to appear, http://arxiv.org/abs/1009.5410.]

\bibitem{App09b}
Appuhamillage, T., V. Bokil, E. Thomann, E. Waymire, B. Wood (2009)
  Solute Transport Across an Interface: A Fickian Theory for Skewness in Breakthrough Curve,
  Water Resour. Res, 46, W07511, doi:10.1029/2009WR008258.


\bibitem{AppShel}
Appuhamillage, T.A., D. Sheldon (2011)
First passage time of skew Brownian motion,
(preprint) http://arxiv.org/pdf/1008.2989

\bibitem{Berkowitz09}
Berkowitz, B., A. Cortis, I. Dror, and H. Scher (2009), Laboratory experiments on dispersive transport across interfaces: The role of flow direction, Water Resour. Res., 45, W02201, doi:10.1029/2008WR007342.



\bibitem{hoteit}
Hoteit, H., R.  Mose, A. Younes, F. Lehmann,
Ph. Ackerer  (2002),
Three-dimensional modeling of mass transfer
in porous media using the mixed hybrid finite
elements and random walk methods,
{\it Mathematical Geology}, {\bf 34}(4), 435-456.

\bibitem{kuo_1999}
Kuo, R. K. H. and Irwin, N. C. and Greenkorn, R. A. and Cushman, J. H. (1999), Experimental investigation of mixing in aperiodic heterogeneous porous media: Comparison with stochastic transport theory, Transport in Porous Media, 37, 169-182.


\bibitem{Lejay08}
Lejay, A., M. Martinez (2008),
A scheme for simulating one-dimensional diffusion
processes with discontinuous coefficients, {\it Ann. Appld.
Probab.} {\bf 16}(1), 107-139. 




\bibitem{Matano}
Matano, Ricardo P.,  Elbio D. Palma. On the upwelling of downwelling currents, 
{\it Journal of Physical 
Oceanography}, {\bf 38}(11):2482 Ð 2500. 


\bibitem{Ramirez06}
Ramirez, J., E. Thomann, E. Waymire, R. Haggerty, B. Wood
(2006), A generalized Taylor-Aris formula and skew diffusion,
{\it SIAM Multiscale Modeling and Simulation} {\bf 5} 3,
786-801.

\vskip 5in

\bibitem{Ramirez08}
Ramirez, J. M., E. A. Thomann, E. C. Waymire, J. Chastanet, and B. D. Wood (2008), A note on the theoretical foundations of particle tracking methods in heterogeneous porous media, Water Resour. Res., 44, W01501, doi:10.1029/2007WR005914. 


\bibitem{Ramirez11a}
Ramirez, J.M. (2011), Multi-skewed
Brownian motion and diffusion in layered media, {\it Proceedings of the American Mathematical Society}. (to appear).

\bibitem{Ramirez11b}
Ramirez, J.M. (2011), Population persistance
under advection-diffusion in river networks,
preprint.


\bibitem{RY}
Revuz, Daniel, and Marc Yor (1999), {\it  Continuous martingales and Brownian
motion} Springer, NY.


\bibitem{Schultz}
Schultz, Cheryl B., and Elizabeth E. Crone (2005), Patch size and connectivity thresholds for butterfly habitat restoration, {\it Conservation
Biology}, {\bf 19}(3), 887-896.

\bibitem{SV}
Stroock, Daniel, W., and S.R. Srinivasa Varadhan (1997), {\it Multidimensional diffusion 
processes} Springer, NY.


\end{small} 
\end{thebibliography}
\end{document}